\documentclass{amsart}
\usepackage{palatino,hhline, bbm}
\usepackage{amsthm, amsmath}
\usepackage{amssymb} 
\usepackage{mathrsfs}
\usepackage[hmargin=3cm, vmargin=3cm]{geometry}
\usepackage[breaklinks,bookmarksopen,bookmarksnumbered]{hyperref}
\usepackage[all]{xy}

\theoremstyle{plain}
\newtheorem*{thm}{Theorem}

\newtheorem*{prop}{Proposition}

\theoremstyle{remark}

\newcommand\pr{\noindent\textit{Proof} : }

\newcommand\rond{\kern 1pt{\scriptstyle\circ}\kern 1pt}

\newcommand\Ker{\operatorname{Ker}}

\newcommand{\mo}{\smallsetminus}

\newcommand\C{\mathbb{C}}
\renewcommand\P{\mathbb{P}}

\renewcommand\O{\mathcal{O}}

\newcommand\iso{\vbox{\hbox to .8cm{\hfill{$\scriptstyle\sim$}\hfill}
\nointerlineskip\hbox to .8cm{{\hfill$\longrightarrow $\hfill}} }}
\newcommand\bir{\vbox{\hbox to .8cm{\hfill{$\scriptstyle\sim$}\hfill}
\nointerlineskip\hbox to .8cm{{\hfill$\dasharrow $\hfill}} }}

\begin{document}
\title{Hyperplane sections of cubic threefolds}
\author[Arnaud Beauville]{Arnaud Beauville}
\address{Universit\'e C\^ote d'Azur\\
CNRS -- Laboratoire J.-A. Dieudonn\'e\\
UMR 7351 du CNRS\\
Parc Valrose\\
F-06108 Nice cedex 2, France}
\email{arnaud.beauville@univ-cotedazur.fr}
 \thanks{I am indebted to A. Dimca for pointing out the reference \cite{A-R}.}
\begin{abstract}
Let $X\subset \P^4$ be a smooth cubic hypersurface. We prove that a general cubic surface is isomorphic to a hyperplane section of $X$. 
\end{abstract}
\maketitle 
\section{Introduction}
A classical result of Sylvester states that a general cubic surface $S\subset \P^3$ admits an equation of the form $L_0^3+\ldots +L_4^3=0$, where the $L_i$ are linear forms (see e.g.\ \cite[\S 9.4]{D}). An equivalent formulation is that $S$ is isomorphic to a hyperplane section of the Fermat cubic threefold $X_0^3+\ldots +X_4^3=0$. We will show in this note that this actually holds for \emph{any} smooth cubic threefold. Surprisingly, the proof reduces to a particular case of the \emph{weak Lefschetz property}  for the Jacobian ring of $X$ --- a case which miraculously   is treated in \cite{A-R}.

\section{Proofs}
Let $X\subset \P^4=\P$ be a smooth cubic hypersurface, defined by an equation $F=0$. 
Let $X^*\subset \P^*$ be the dual hypersurface. For any $H\in\P^*\mo X^* $, the hyperplane section $X\cap H$ is a smooth cubic surface; we get in this way a morphism $s^{}_X$ from $\P^*\mo X^* $ into the moduli space $\mathscr{M}_3$ of smooth cubic surfaces. Our aim is to determine when this map is dominant. 

Let $R$ denote the graded ring $\C[X_0,\ldots ,X_4]$.  Let
$J$ be the ideal $(F'_0,\ldots ,F'_{4})$ in  $R$, where we put $F'_i:=\dfrac{\partial F}{\partial X_i} $, and let
$\mathfrak{J}:= R/J$ be the Jacobian ring of $F$.

\begin{prop}
Let $H\in \P^*\mo X^*$, given by a linear form $L\in R_1$. The map $s^{}_X: \P^*\mo X^* \rightarrow \mathscr{M}_3$ is \'etale at $H$ if and only if the multiplication map $\times L: \mathfrak{J}_2\rightarrow \mathfrak{J}_3$ is injective. 
\end{prop}
\pr We put  $S:=X\cap H$. The tangent space to $\P^*$ at $H$ is naturally isomorphic to $H^0(H,\O_{H}(1))=\allowbreak H^0(S,\O_S(1))$, and the tangent map to $ s^{}_X$ at $H$ is the coboundary map
\[\partial : H^0(S,\O_S(1))\rightarrow H^1(S,T_S)\]
deduced from the exact sequence
\[0\rightarrow T_S\rightarrow T_{X|S}\rightarrow \O_S(1)\rightarrow 0\,.\]
Since $\dim H^0(S,\O_S(1)) = \dim H^1(S,T_S)= 4$ and $H^0(S,T_S)=0$, $s^{}_X $ is \'etale at $H$ if and only if $H^0(S, T_{X|S})=0$.

By restricting to $S$ the exact sequence
\[0\rightarrow T_X\rightarrow T_{\P|X}\rightarrow \O_X(3)\rightarrow 0\,,\]
we see that $H^0(S, T_{X|S})$ is the kernel of the map $\varphi : H^0(S, T_{\P|S})\rightarrow H^0(S,\O_S(3))$. We choose the coordinates in $\P$ so that $L=X_0$. 
The space $H^0(S, T_{X|S})$ is generated by the vector fields (restricted to $S$) $\ X_i\dfrac{\partial }{\partial X_j} \ (i>0, j\geq 0)$,  with the relation $\sum\limits_{i\geq 1} X_i\dfrac{\partial }{\partial X_i}=0 $; we have $\varphi (X_i\dfrac{\partial }{\partial X_j} )=X_i F'_j$. 

Let  $V=\sum\limits_{i=0}^4 L_i  \dfrac{\partial }{\partial X_i}$ be  a nonzero element of $H^0(S, T_{X|S})$, where the $L_i$  are linear forms in $X_1,\ldots ,X_4$. Assume that $\varphi (V)=0$. 
This means that there exists $Q\in R_2$ and $a\in\C$ such that  
\begin{equation}\label{1}
\sum_{i=0}^4 L_i F'_i= X_0Q+3aF\,;
\end{equation}
this implies that the class of $Q$ in $\mathfrak{J}_2$ is mapped to zero in $\mathfrak{J}_3$ under multiplication by $X_0$.

Let us show that this class is nonzero.
Suppose  $Q\in J_2$, so that $Q=\sum a_iF'_i$. Using the Euler relation $\sum X_iF'_i=3F$, (\ref{1}) becomes 
\[\sum_{i=0}^4 (L_i-a_iX_0-aX_i)F'_i=0\,.\]
Now the only nontrivial relation $\sum M_iF'_i=0$ with the $M_i$ in $R_1$ is the Euler relation
 (this follows from the cohomology exact sequence associated to   
$0\rightarrow  T_X\rightarrow T_{\P|X}\rightarrow \O_X(3)\rightarrow 0$, plus the fact that $H^0(X,T_X)=0$). 

Therefore there exists $b\in\C$ such that $L_i-a_iX_0-aX_i=bX_i$ for all $i$; this implies $a_i=0$ for $i\geq 1$, hence $L_i= (a+b)X_i$, and $L_0=0$. Thus $V= (a+b)\sum\limits_{i\geq 1} X_i\dfrac{\partial }{\partial X_i} =0$, a contradiction.

\medskip	
Conversely, if $\times X_0: \mathfrak{J}_2\rightarrow \mathfrak{J}_3 $ is not injective, there exist forms $Q\in R_2\mo J_2$ and $L_0,\ldots ,L_4$ in $R_1$ such that 
\begin{equation}\label{2}
X_0Q=\sum L_iF'_i\,.
\end{equation}
Replacing $Q$ by $Q-\sum a_iF'_i$, where $a_i$ is the coefficient of $X_0$ in $L_i$, 
we can assume that the $L_i$ are linear forms in $X_1,\ldots ,X_4$. Then the vector field $V:= \sum L_i\dfrac{\partial }{\partial X_i}$ satisfies $\varphi (V)=0$. 

We just have to check that $V$ is nonzero. If it is, there exists $b\in\C$ such that $L_i=bX_i$ for $i\geq 1$ and $L_0=0$. Then (\ref{2}) becomes 
\[ X_0(Q+bF'_0)= \sum_{i\geq 0}bX_iF'_i=3bF\,.\]
Since $F$ is irreducible, this implies $b=0$ and $Q=0$, contradiction.\qed

\medskip	
\noindent\emph{Remark}$.-$
It is  easy to give examples of pairs $(X,H)$ such that $s^{}_X$ is not \'etale at $H$. 
For instance if $X$ is the Fermat cubic and $H$ is the hyperplane $X_0=0$, the elements $X_i\dfrac{\partial }{\partial X_0}\ (i\geq 1)$ of $H^0(S, T_{\P|S})$ belong to $\Ker \varphi $, so the tangent map to $s_X$ at $H$ is zero. In fact, $s^{}_X$ contracts the lines in $\P^*$ consisting of the hyperplanes $X_0=tX_i$, $t\in\C$, $i\geq 1$.

\begin{thm}
The map $s^{}_X: \P^*\mo X^* \rightarrow \mathscr{M}_3$ is dominant.
\end{thm}
\pr By \cite[Theorem 1]{A-R}, the multiplication map $\times L: \mathfrak{J}_2\rightarrow \mathfrak{J}_3$ is injective for $L$ general in $R_1$
(weak Lefschetz property in degree 2). Thus  by the Proposition,   $s^{}_X$ is \'etale at a general point of $\P^*$, hence dominant.\qed

\medskip	
\noindent\emph{Remarks}$.-\ 1)$
The proof of the Proposition applies identically to  a smooth hypersurface $X$ of degree $d$ in $\P^{n}$, with $d,n\geq 3$. With the previous notation, we get that \emph{$s^{}_X$ is unramified at $H$ if and only if the multiplication map $\times L:\mathfrak{J}_{d-1}\rightarrow \mathfrak{J}_d$ is injective}.  If $\mathfrak{J}$ satisfies the weak Lefschetz property in degree $d-1$, we  conclude that $s^{}_X$ \emph{ is generically finite onto its image}. This has been proved recently in  \cite{B-M} for $d\geq n+2$.
The cases $3<d<n+2$ and $d=3$, $n\geq 5$ remain open.

\smallskip	
2) The proof of the Proposition works over any algebraically closed field $k$. However, the weak Lefschetz property used in the Theorem does not always hold. For instance if 
$X$ is the Fermat cubic threefold in characteristic 2,  the Jacobian ideal is generated by the squares of all linear forms; if $\ell$ and $m$ are linearly independent linear forms, we have $\ell\cdot m\neq 0$ in $\mathfrak{J}_2$, but $\ell\cdot (\ell\cdot m)=0$ in $\mathfrak{J}_3$, so the map $\times \ell:\mathfrak{J}_2\rightarrow \mathfrak{J}_3$ is not injective.

In fact in this case all smooth hyperplane sections of $X$ are isomorphic, see \cite[Th\'eor\`eme 2]{B}.

\end{document}